\documentclass[11pt]{amsart}
\usepackage{geometry} 
 \usepackage[english]{babel}
\usepackage{amsmath, amscd, amsthm, amssymb, mathrsfs}
\usepackage{verbatim,alltt}
\usepackage{color}
\usepackage{caption}
\usepackage{subcaption}
\usepackage{lscape,supertabular}
\usepackage{tikz}

\theoremstyle{plain}

\theoremstyle{definition}

\newcommand{\Mo}{\mathrm{M_{11}}}
\newcommand{\Md}{\mathrm{M_{12}}}
\newcommand{\Mvd}{\mathrm{M_{22}}}
\newcommand{\Mvt}{\mathrm{M_{23}}}
\newcommand{\Mvq}{\mathrm{M_{24}}}
\newcommand{\McL}{\mathrm{McL}}
\newcommand{\He}{\mathrm{He}}

\newcommand{\Cod}{\mathrm{Co_2}}
\newcommand{\Cot}{\mathrm{Co_3}}
\newcommand{\Fivd}{\mathrm{Fi_{22}}}

\newcommand{\Suz}{\mathrm{Suz}}
\newcommand{\Ru}{\mathrm{Ru}}
\newcommand{\Ju}{\mathrm{J_1}}
\newcommand{\Jd}{\mathrm{J_2}}
\newcommand{\Jt}{\mathrm{J_3}}
\newcommand{\HS}{\mathrm{HS}}

\newcommand{\ON}{\mathrm{O'N}}

\title{An Atlas of subgroup lattices of finite almost simple groups}
\author{Thomas Connor}
\address{Universit\'e Libre de Bruxelles,
D\'epartement de Math\'ematiques - C.P.216,
Boulevard du Triomphe,
B-1050 Bruxelles, Boursier FRIA}
\email{tconnor@ulb.ac.be}

\author{Dimitri Leemans}
\address{University of Auckland,
Department of Mathematics,
Private Bag 92019,
Auckland, New Zealand}
\email{d.leemans@auckland.ac.nz}

\date{\today} 
\begin{document}
\begin{abstract}
We announce the publication of an atlas of subgroup lattices for a large collection of finite almost simple groups. This atlas is made available online.
\end{abstract}

\maketitle
{\bf Keywords:} Simple group, permutation group, subgroup lattice, computer algebra, taxonomy.

\section{Introduction}
The Classification of the Finite Simple Groups (CFSG) emphasizes the importance of the finite simple groups in Group Theory. It is one of the most impressive achievements in the history of Mathematics. We refer to~\cite{Wilson09} and the references provided there for a broad literature on this wonderful theorem.
Among the amazing achievements in this branch of Mathematics, we find the {\sc Atlas} of Finite Groups~\cite{Atlas} as well as the online version of the {\sc Atlas} of Finite Group Representations~\cite{onlineAtlas}.

Over the years, the finite simple groups have received a lot of attention with respect to the study of geometry. The Theory of Buildings due to Jacques Tits, who was awarded the Abel Prize in 2008, illustrates this perfectly. We refer for instance to~\cite{Abramenko08} and references provided there.
Much work has been done in this respect with the study of incidence geometries associated to finite almost simple groups (we refer to~\cite{Buek95,BuekCohen,DL08} and references cited there for a large documentation on this aspect).

The computations of subgroup lattices and tables of marks of permutation groups, and in particular sporadic simple groups, have been a subject of interest for many decades, linked among others to the search for a unified geometric interpretation of all finite simple groups. 
{
Joachim Neub\"user gave in~\cite{Neu1960} the first algorithm that was implemented later on in the computational software {\sc Cayley} and its successor {\sc Magma}.
Francis Buekenhout computed in 1984 the lattices of $\Mo$ and $\Ju$~\cite{Bue86}. Then, Herbert Pahlings did the lattice of $\Jd$~\cite{Pah87} in 1987. In 1988, Buekenhout and Sarah Rees produced the lattice of $\Md$ (see~\cite{BR88} and~\cite{DL99} for a few corrections). In 1991, Pfeiffer computed the table of marks of  $\Jt$ and in 1997, those of $\Mvd$, $\Mvt$ $\Mvq$,  $\McL$~\cite{Pfe97}. Also in 1997, Merkwitz got the tables of marks of  $\He$ and  $\Cot$.
In 1998, Derek Holt computed all conjugacy classes of subgroups of $\ON$ (personal communication).
In a more general setting again, John Cannon, Bruce Cox and Derek Holt described in~\cite{CCH2001} a new algorithm to compute the conjugacy classes of subgroups of a given group that was used in {\sc Magma} until 2005.
Progresses on the computation of maximal subgroups of a given group by Cannon and Holt~\cite{CH2004} led Leemans to a much faster algorithm to compute the subgroup lattice of a given group that is now available in {\sc Magma}.
In 2007, Leemans computed the full subgroup lattices of $\HS$,  $\Ru$,  $\Suz$,  $\ON$,  $\Cod$ and  $\Fivd$ using permutation degree reduction at each step of the computation. 
}
Recently, Naughton and Pfeiffer produced a new algorithm to compute the table of marks of a cyclic extension of a group~\cite{NP2012}.
 
The knowledge of the subgroup lattice of a group $G$ is a powerful tool to study the symmetrical objects on which $G$ acts. For instance, in \cite{Buek13}, \cite{Connor:2013} and ~\cite{Lee2010}, the authors build flag-transitive coset geometries of ranks 2, 3 and 5 for $\ON$ by identifying boolean lattices in the subgroup lattice of $\ON$. In~\cite{ConnorLeemans2013}, the authors develop an algorithm in order to count the number of regular maps on which a finite group $G$ acts regularly. This algorithm makes an intensive use of the knowledge of the subgroup lattices of $G$. Then they illustrate the algorithm on the group $\ON$. In a more general approach, \cite{Downs:1987} discusses the problem of enumerating regular objects with a given automorphism group. The authors introduce, among other things, the M\"obius function for a group $G$ as a tool to enumerate regular objects. The knowledge of the M\"obius function of $G$ relies on the knowledge of the full subgroup lattice of $G$. 

In the spirit of contributing to the study of the finite simple groups, we present here an algorithm that determines the subgroup lattice of a given permutation group. This algorithm, designed by Leemans in 2007, proves itself to be lighter in memory and sometimes even faster than the already implemented function {\tt SubgroupLattice} in the software {\sc Magma}~\cite{Magma}. We also present an algorithm to determine structures of groups with a computer. Those two algorithms allow us to produce an atlas of subgroup lattices for a large number of finite almost simple groups. The atlas is made available online at

\begin{center}
{\sf http://homepages.ulb.ac.be/$\sim$tconnor/atlaslat.}
\end{center}

Our paper is organized as follows. In section~\ref{slat}, we present {the} algorithm that computes the subgroup lattice of a permutation group $G$. It features a systematic reduction of the permutation degree of the subgroups of $G$ and a possibility to start the computation of the subgroup lattice of $G$ with partial information of the lattice already available. 
This permits to compute the subgroup lattice of very large groups, like the O'Nan sporadic group $\ON$ or even its automorphism group, currently out of reach with the {\tt SubgroupLattice} function of {\sc Magma}, but also $\Cod$ and $\Fivd$.
In section~\ref{struct} we discuss the problem of describing the structure of a group in an efficient way. We present an algorithm that provides a structure for a group, based on a choice of suitable normal subgroups. 
In section~\ref{results} we present our atlas of subgroup lattices as an application of the algorithms presented and discussed in this paper.

\section{The subgroup lattice of a permutation group}\label{slat}

We refer to~\cite{Hall:1976} as a reference on subgroup lattices. A lattice is a partially ordered set, or poset, any two of whose elements $a,b$ have a least upper bound $a\cup b$ and a greatest lower bound $a\cap b$. The subgroups of a group $G$ may be taken as the elements of a lattice $L(G)$ under the operations of union and intersection. The poset of conjugacy classes of subgroups forms also a lattice: two conjugacy classes $A$ and $B$ are such that $A \supseteq B$ provided that any subgroup of $B$ is contained in some subgroup of $A$. We call this lattice the subgroup lattice of $G$, rather than the lattice of conjugacy classes of subgroups of $G$ for the sake of brevity, and we denote it with $\Lambda(G)$. Our terminology is also the one used in {\sc Magma}. This lattice can be refined with the length of each conjugacy class of subgroups. Moreover, given two conjugacy classes of subgroups $A\supset B$, we define $n_{AB}$ to be the number of subgroups of class $B$ contained in any subgroup of class $A$; alike we define $n_{BA}$ to be the number of subgroups of class $A$ containing a subgroup of class $B$. Consider the set $N$ of numbers $n_{XY}$ for every couple of classes $\{X,Y\}$ such that $X\subset Y$ or $X \supset Y$ and there does not exist $Z$ such that $X\subset Z \subset Y$ or $X \supset Z \supset Y$. The subgroup lattice $\Lambda(G)$ together with the length of each conjugacy class and the set $N$ is called the weighted subgroup lattice of $G$.

We describe in this section a powerful and natural algorithm to compute the weighted subgroup lattice of a given group $G$. The correctness of this algorithm is obvious.

Start with a set {\tt classes} which is empty and a set {\tt sgr} containing just one element, 
namely the group $G$ for which we want to compute the subgroup lattice.
While {\tt sgr} is nonempty, pick one element $H$ out of {\tt sgr} and put it in {\tt classes}. Obviously, it is $G$ the first time.
Reduce the permutation degree of $H$ and let $\phi : H\rightarrow \tilde{H}$ be an isomorphism between $H$ and $\tilde{H}$ where $\tilde{H}$ has a reduced permutation degree.
Compute the maximal subgroups of $\tilde{H}$ and for each maximal $\tilde{M}$, add $M := \phi^{-1}(\tilde{M})$ to {\tt sgr} provided there is no subgroup in {\tt sgr} conjugate to $M$ in $G$.
During that process, keep track of inclusions of respective subgroups considered.
At the end of this process, in {\tt classes} there is one representative of each conjugacy class of subgroups of $G$. 
Moreover, we also have the maximal inclusions between classes. 
So the subgroup lattice is determined.
The weighted subgroup lattice can be determined in the process by computing weighted inclusions at each step.

A {\sc Magma} implementation of the algorithm described above to compute the subgroup lattice of a given group is available on the webpage of the atlas. Observe that we use the {\tt DegreeReduction} function in {\sc Magma} to get $\phi$ and $\tilde{H}$ for every subgroup $H$ above. 
This improvement can save a lot of time and memory.
For instance, consider $L_3(7):2$, one of the maximal subgroups of the O'Nan sporadic group $\ON$, acting on 122760 points (the smallest permutation representation of $O'N$). Then {\sc Magma} v.2.19 needs 13 seconds and more than 200 Mb of memory to compute its maximal subgroups on a computer running at 2.9~GHz. If we reduce the degree of $L_3(7):2$ on 5586 points by using the {\tt DegreeReduction} function, then {\sc Magma} computes them in less than half a second and takes about 20 Mb of memory.

\section{The structures of a group}\label{struct}

\subsection{Preliminary remarks}

Given any finite group $G$, it is always desirable to identify $G$ in some sense. This identification can be done for instance in a geometrical way by determining the action of $G$ on some set or by algebraic means. In particular, most finite simple groups can be named after their action on some structured set or after the mathematician that discovered them (like the Suzuki groups or most of the sporadic groups). However some groups carry very different names, depending on the incarnation of the group that the context requires to emphasize. This is the case for instance of $U_4(2)$. Indeed, $$S_{4}(3) \cong U_4(2) \cong O_5(3) \cong O^-_6(3) \cong W(E_6).$$ Each of the names of this group emphasizes one of its actions on a structured set of particular interest. Therefore, when speaking about this group, one has to choose carefully the name that should be used depending on the context. This observation means that one has to be aware of possible isomorphisms between different incarnations of a group.

In {\sc Magma}, there exists a database of finite simple groups. Given a simple group $G$, one can thus ask {\sc Magma} to name $G$ by using the function {\tt NameSimple}. This function returns a triple of integers that permits to identify $G$ as a group of one of the infinite families of finite simple groups, or as one of the sporadic groups.
Many non simple groups can also be identified in a canonical way. This is the case of most of the almost simple groups for instance, but also the case of the dihedral groups, or the groups $AGL(n,q)$. Abelian groups are also identified easily by a name thanks to the classification theorem of abelian groups.
However, most of the finite groups are not almost simple, and identifying them in an efficient way by a name can be tricky. For instance, Leemans exhibited two non isomorphic primitive groups in~\cite{Leemans97} that satisfy the following property: they have isomorphic posets of conjugacy classes of subgroups and for each normal subgroup $N$ of the first, there is a normal subgroup isomorphic to $N$ in the second group such that the quotients by $N$ are isomorphic. In other words, it is not possible to make a difference between those two groups by giving them names based on any quotient by a normal subgroup. This shows that the taxonomy of groups is a difficult and possibly not solvable problem. Hence we should not look for a deterministic algorithm that gives names to groups since it is readily impossible.

The case of $p$-groups is also particularly difficult to handle. For instance, there are roughly 50 billions pairwise non-isomorphic groups of order 1024 and hence, finding a way to give distinct names to each of them is hopeless, unless we decide to assign a number to each of them, as is done for instance in the {\sc SmallGroups} database provided by~\cite{OBrien01}.

\subsection{Algorithmic approach}
Let $G$ be a group and let $N$ be a normal subgroup of $G$. Denote by $Q$ the quotient group $G/N$. Then $G$ can be written as $N.Q$ where the dot \mbox{`` . "} denotes an extension that can be split (that is, a direct or a semi direct product) or non split. We denote a direct product by `` $\times$ ", a semi direct product by `` $:$ " and a non split extension by `` $\cdot$ ".

We recall that a {\em composition series} for $G$ is a sequence of subgroups $H_i$, $i\in \{0,\ldots,n+1\}$ such that $$1=H_0 \triangleleft H_1 \triangleleft H_2 \triangleleft \ldots \triangleleft H_{n} \triangleleft H_{n+1} = G$$ where all inclusions are strict, i.e. $H_i$ is a maximal normal subgroup of $H_{i+1}$. This is equivalent to require that the {\em composition factors} $Q_i=H_{i+1}/H_{i}$ are simple groups, $i=0, \ldots, n$. Clearly the group $G$ can be written $H_{n}.Q_n$. Alike, $H_n$ can be written $H_{n-1}.Q_{n-1}$ and thus $G$ can be written $(H_{n-1}.Q_{n-1}).Q_n$. Proceeding inductively, we can finally write $$G=(\ldots(Q_0.Q_1).Q_2)\ldots).Q_n.$$ However in order to reduce the notations, we always suppose that the products are left associate and we can thus avoid to write parentheses whenever there is no possible confusion. Therefore by $G\cong A.B.C$ we mean $G\cong(A.B).C$.

The Jordan--H\"older theorem states that every finite group has a unique composition series up to the order of the terms~\cite{Wilson09}. Obviously, two non isomorphic groups can have the same composition series. This is the case for instance of $S_5$ and $A_5 \times 2$. In this particular example, it is not enough to use the composition series of those two groups to distinguish them. However $S_5\cong A_5:2$ but $S_5 \ncong A_5\times 2$.

On basis of the previous observations, we detail an algorithm that produces a name for a group $G$ in terms of a product of its composition factors. We detail afterwards an improved algorithm that we actually used in order to produce the lattices of our atlas.
First of all, given $N\triangleleft G$ we need to check whether the extension $N.Q$ is split or not, i.e. $G\cong N:Q$ or $G\cong N\cdot Q$, where $Q \cong G/N$ as usual. 
%
If $N$ is a maximal normal subgroup of $G$, then $Q$ is simple. We can use the database of simple groups in {\sc Magma} to identify $Q$ and give it a name. We can now easily extract the following algorithm from the previous observations. If $G$ is simple, we are done. Suppose $G$ is not simple. Compute a composition series of $G$ and the corresponding composition factors. At step $n-i+1$, identify the simple group $Q_i$ and check whether $H_i.Q_i$ is split or non split.  If it is split, check moreover if the extension is a direct product. The procedure returns the group $G$ written as $G\cong Q_0._0Q_1._1\ldots Q_{n-1}._{n-1}Q_n.$ where $._i$ is a symbol in  $\{\times, :, \cdot\}$. Applying this procedure to $S_5$ for instance would produce $A_5:2$. Applying it to the dihedral group $D_{40}$ would produce $5\times 2\cdot 2:2$. Unfortunately, this could also be the result after applying this algorithm to $5 \times D_8$. Finally applying it to an elementary abelian group $2^5$ would produce $2\times 2\times 2 \times 2\times 2$.

There is an obvious improvement of this algorithm. The guideline is that some nonsimple groups can be identified in a canonical way like the symmetric groups or the dihedral groups for instance. Moreover in the process of building the Atlas that we describe in this article, we observed for example that the group $S_3\times S_3$ would not be named correctly most of the time, or $A_4$ would be written $2^2:3$. Therefore we produced a database of selected groups that our algorithm checks prior to computing the list of normal subgroups of $G$. The algorithm also checks possible isomorphisms of $G$ with `classical' groups (like the symmetric groups or the dihedral groups, for instance).
If $G$ is not immediately identifiable, our algorithm computes the list of normal subgroups of $G$. If a normal subgroup $N$ or a quotient $G/N$ is appealing then our algorithm would select it and proceed inductively.

\begin{figure}
\scriptsize
\begin{alltt}
if \(G\) is simple then identify \(G\)
else if \(G\) is in the database then identify \(G\)
else if \(G\) has a `desirable property' then identify \(G\)
else
 compute the list \(L\) of normal subgroups of \(G\)
 for each subgroup \(N\) in \(L\) in decreasing order do
  if \(N\) is simple or has a `desirable property' or is in the database then
   identify \(N\) and identify the extension between \(N\) and \(G/N\)
   proceed inductively on \(G/N\)
  else if \(G/N\) has a `desirable property' or is in the database then
   identify \(G/N\) and identify the extension between \(N\) and \(G/N\)
   proceed inductively on \(N\)
  else if no \(N\) and no \(G/N\) is desirable then
   take the largest \(N\) and proceed inductively
\end{alltt}
\caption{An improved naming algorithm}
\label{betteralgo}
\end{figure}
\section{The atlas}~\label{results}
For every almost simple group of order at most 1,000,000 appearing in the online version of the Atlas of Finite Groups~\cite{onlineAtlas}, we computed its subgroup lattice with the {\sc Magma} implementation of our algorithm, available on the homepage mentioned below. Given such a lattice $\Lambda$, we ran the algorithm described in Figure~\ref{betteralgo} on every subgroup in $\Lambda$. We also proceeded in this way for some groups of order larger than 1,000,000 like some large sporadic groups. The result is an atlas of more than a hundred subgroup lattices of almost simple groups with a structure provided for every subgroup of each group.
The atlas of subgroup lattices is available online at

\begin{center}
{\sf http://homepages.ulb.ac.be/$\sim$tconnor/atlaslat.}
\end{center}

\noindent Groups are subdivided in several families, namely almost simple groups of sporadic type, alternating type, linear type, symplectic type, orthogonal type, unitary type and exceptional Lie type.
For each group $G$ in the atlas, a {\tt pdf} file containing the subgroup lattice of $G$ is available for download.

\section{Acknowledgements}

Part of this research was done while the first author was visiting the second author at the University of Auckland. He therefore gratefully acknowledges the University of Auckland for its hospitality. He also acknowledges the {\em Fonds pour la formation \`a la Recherche dans l'Industrie et l'Agriculture}, (F.R.I.A.) and the {\em Fonds National pour la Recherche Scientifique} (F.N.R.S), national Belgium science fundings, for financial support.
The second author acknowledges financial support of the Royal Society of New Zealand Marsden Fund (Grant 12-UOA-083).

\bibliographystyle{plain}

\end{document}